\documentclass[12pt]{article} 
\begin{document}

%Self-Definitions  
\def\endproof {\rule[-0.5mm]{1.2ex}{1.2ex}}
\def\C{{\bf C}}
\def\ord{{\rm ord} }
\def\mult{{\rm mult}}

\author{Nguyen Van Chau \thanks{
Supported in part by the National Basic Program of Natural Science, Vietnam. }}

\date{}
\title{ A Simple Proof of Jung' Theorem  on Polynomial Automorphisms of $\C^2$}
\maketitle

{\small {\bf Abstract.} The Automorphism Theorem, discovered first by Jung in 1942, asserts that if $k$ is a field, then every  polynomial automorphism of $k^2$ is a finite product of linear automorphisms and automorphisms of the form $(x,y)\mapsto(x+p(y), y) $ for $p\in k[y]$. We present here a simple proof for the case $k=\C$ by using Newton-Puiseux expansions.}

\bigskip
\noindent {\bf 1.}
In this note we present a simple proof of the following theorem on the structure of the group $GA(\C^2)$ of polynomial automorphisms of $\C^2$

\medskip

\noindent{\bf Automorphism Theorem. }{\it Every polynomial automorphism of $\C^2$ is tame, i.e. it is a finite product of linear automorphisms and automorphisms of the form $(x,y)\mapsto (x+p(y),y)$  for one-variable polynomials $p\in \C [y]$.}

\medskip

This theorem was first discovered by Jung [J] in 1942. In 1953, Van der Kulk [Ku] extended it to a field of arbitrary characteristic. In an attempt to understand the structure of $GA(\C^n)$ for large $n$, several proofs of Jung's Theorem have presented by Gurwith [G], Shafarevich [Sh],  Rentchler [R], Nagata [N], Abhyankar and Moh [AM], Dicks [D], Chadzy'nski and  Krasi'nski [CK]  and  McKay and Wang [MW] in different approaches. They are related to the mysterious Jacobian conjecture, which asserts that a polynomial map of $\C^n$ with non-zero constant Jacobian is an automorphism. This conjecture dated back to 1939 [K], but it is still open even for $n=2$. We refer to [BCW] and [E] for nice surveys on this conjecture. 

\medskip
\noindent{\bf 2.} The following essential observation due to van der Kulk [Ku] is the crucial step in some proofs of Jung' theorem.

\medskip

\noindent{\bf Division Lemma: } {\it $F=(P,Q) \in GA(\C^2) \Rightarrow \deg P\vert \deg Q \ or  \deg Q\vert \deg P.$}

\medskip
Abhyankar and Moh in [AM] deduced it as a consequence of the theorem on the embedding of a line to the complex plane. McKay and Wang [MW] proved it by using formal Laurent series and the inversion formula. Chadzy'nski and  Krasi'nski in [CK] obtained the Division Lemma from a formula of geometric degree of polynomial maps $(f,g)$ that the curves $f=0$ and $g=0$ have only one branch at infinity. Here, we will prove this lemma by examining the intersection of irreducible branches at infinity of the curves $P=0$ and $Q=0$ in term of Newton-Puiseux expansions. 

Our proof presented here is quite elementary and simpler than any proof mentioned above. It uses the following  two elementary facts on Newton-Puiseux expansions (see, for example,  [BK]). 

Let $h(x,y)=y^n+a_1(x)y^{n-1}+\dots +a_n(x)$ be a reducible polynomial. Looking in the compactification $\C P^2$ of $\C^2$, the curve $h=0$ has some irreducible branches located at some points in the line at infinity, which are called the {\it irreducible branchs at infinity}. For  such a branch $\gamma$, the Newton' algorithm allows us to find a meromorphic parameterization of $ \gamma$, an one-to-one meromorphic map $ t \longmapsto (t^m,u(t)) \in \gamma $  defined for $t$ large enough,
$$ 
u(t)=t^m\sum\limits_{k=0}^\infty b_kat^{-k},  \   \gcd \{k :b_k \neq 0 \}=1,  
$$ 
The fractional power series $u(x^{\frac 1m})$ is called a {\it Newton-Puiseux expansion at infinity} of $\gamma$ and the natural number $\mult(u):=m$ - the {\it multiplicity } of $u$. 

The first fact is a simple case of Newton's theorem (see in [A]).

\medskip
\noindent{\bf Fact 1. }
 {\it Suppose the curve $h=0$ has only one irreducible branch at infinity  and $u$ is a Newton-Puiseux expansion at infinity of this branch. Then  $$h(x,y)=\prod_{i=1}^{\deg h}(y-u(\epsilon^i x^{\frac1{\deg h}}))$$
and $\mult(u)=\deg h $,
where $\epsilon$ is a primitive $\deg h$-th root of $1$.}

\medskip

Let $\varphi(x,\xi)$ be a  finite fractional power series of the form 
$$\varphi (x,\xi)=\sum_{k=0}^{n_\varphi-1}c_kx^{1-\frac k{m_\varphi }}+\xi x^{1-\frac {n_\varphi}{m_\varphi}},\eqno (1)$$
where $\xi$ is a parameter and $\gcd (\{k=0,\dots n_\varphi-1 : c_k \neq 0 \}\cup \{ n_\varphi\})=1$. Let us represent
$$h(x,\varphi (x,\xi))=x^{\frac{a_\varphi}{m_\varphi}}(h_0(\xi)+\mbox{lower terms in } x^{\frac 1{m_\varphi}}), \ h_0(\xi)\neq 0. \eqno (2)$$

The second fact is deduced from the Implicit Function Theorem. 

\medskip
\noindent{\bf Fact 2. } {\it Let $\varphi$ and $h_0$ be as in (1) and (2). If $c$ is a simple zero of $h_0(\xi)$, then there is a Newton-Puiseux expansion at infinity 
$$u(x^{\frac 1{m_\varphi}})=\varphi (x,c+\mbox{\rm lower terms in } x^{\frac 1{m_\varphi}} )$$
for which $h(x,u(x^{\frac 1{m_\varphi}}))\equiv 0$. Furthermore, $\mult(u)$ divides $m_\varphi$ and $\mult(u)=m_\varphi$ if $c \neq 0$.}

\medskip

\noindent {\bf 3. Proof of the Division Lemma.} Given $F=(P,Q) \in GA(\C^2)$. We may assume that $\deg P > \deg Q$ and  we will prove that $\deg Q$ divides $\deg P$.
By choosing a suitable linear coordinate, we can express  
$$P(x,y)= y^{\deg P}+ \mbox{lower terms in } y$$ 
$$Q(x,y)=y^{\deg Q}+\mbox{lower terms in  } y.$$
Observe that $F$ is a polynomial diffeomorphism of $\C^2$ and 
$$J(P,Q):=P_xQ_y-P_yQ_x \equiv const. \neq 0.$$
Then, $P$ and $Q$ are reducible and each of the curves $P=0$  and $Q=0$ is diffeomorphic to $\C$ which has only one irreducible branch at infinity. Let $\alpha$ and $\beta$ be the unique irreducible branches at infinity of $P=0$ and $Q=0$ , respectively. Then, by Fact 1 we can find Newton-Puiseux expansion $u(x^{\frac 1{\deg P}})$ and $v(x^{\frac 1{\deg Q}})$  with $\mult(u)=\deg P$ and $\mult(v)=\deg Q$ such that
$$P(x,y)=\prod_{i=1}^{\deg P} (y-u(\sigma^ix^{\frac 1{\deg P}}))$$
$$Q(x,y)=\prod_{j=1}^{\deg Q} (y-v(\delta^jx^{\frac 1{\deg P}})),$$ 
where $\sigma$ and $\delta$ are  primitive $\deg P$-th and $\deg Q$-th roots of $1$, respectively.

Put $\theta:=\min_{ij} \ord(u(\sigma^ix^{\frac 1{\deg P}})-v(\delta^jx^{\frac 1{\deg Q}}))$.
Without loss of generality, we can assume 
$\ord(u(x^{\frac 1{\deg P}})-v(x^{\frac 1{\deg Q}}))=\theta$.
We define a fractional power series $\varphi (x,\xi)$ with parameter $\xi$  by deleting in $u$ all terms of order no large than $\theta$ and adding to it the term $\xi x^\theta$,
$$\varphi (x,\xi)=\sum_{k=0}^{n_\varphi-1}c_kx^{1-\frac k{m_\varphi}}+\xi x^{1-\frac {n_\varphi}{m_\varphi}}$$ 
with $\gcd \{k=0,\dots K-1 : c_k \neq 0 \}\cup \{ n_\varphi \}=1,$ 
where $1-\frac {n_\varphi}{m_\varphi}=\theta$. Then, by definition
$$u(x^{\frac 1{\deg P}})=\varphi (x,\xi_u(x)) \mbox{ with } \xi_u(x)=\alpha_u+\mbox{ lower terms in } x,$$
$$v(x^{\frac 1{\deg Q}})=\varphi (x,\xi_v(x)) \mbox{ with }  \xi_v(x)=\beta_v+\mbox{ lower terms in } x$$
and $\alpha_u-\beta_v \neq 0.$
Let us  represent  
$$ 
P(x,\varphi (x,\xi ))=x^{\frac {a_\varphi}{m_\varphi} }(P_\varphi (\xi )+\mbox{ lower terms in } x^{\frac 1{m_\varphi}  })
$$ 
$$ 
Q(x,\varphi (x,\xi ))=x^{\frac {b_\varphi}{m_\varphi }}(Q_\varphi (\xi )+\mbox{ lower terms in } x^{\frac 1{m_\varphi}  })
$$
where $a_\varphi$ and $b_\varphi$ are integers and $0\neq P_\varphi , Q_\varphi \in \C[\xi]$.

\medskip

\noindent{\bf Claim 1. }\\
(a) {\it $P_\varphi(\alpha_u)=0$ and $Q_\varphi(\beta_v)=0$.}\\ 
(b) {\it  The polynomials  $P_\varphi(\xi)$ and $Q_\varphi(\xi)$ have no common zero.}

\medskip
\noindent{\it Proof.} (a) is implied from the equalities $P(x,\varphi (x,\xi_u(x) ))=0$ and \\ $Q(x,\varphi (x,\xi_v(x) ))=0$. For (b), if $P_\varphi(\xi)$ and $Q_\varphi(\xi)$ have a common zero $c$, then by Fact 2 there exists series 
$$\bar \xi_u(x)=c+\mbox{ lower terms in } x,$$ 
$$ \bar\xi_v(x)=c+\mbox{ lower terms in } x$$
such that $\varphi (x,\bar\xi_u(x) )$ and $\varphi (x, \bar\xi_u(x) )$ are Newton-Puiseux expansions at infinity of $\alpha$ and $\beta$, respectively. For these expansions $\ord (\varphi (x,\bar \xi_u(x) )-\varphi (x, \bar\xi_v(x) )<\theta$. This contradicts to the definition of $u$ and $v$. \endproof
\medskip

\noindent{\bf Claim 2. }{\it $P_\varphi$ and $Q_\varphi$ have only simple zeros.}

\medskip
\noindent{\it Proof.}
First, observe that
$$a_\varphi >0, \ b_\varphi> 0.\eqno(3)$$
Indeed, for instance, if  $a_\varphi \leq 0 $, then $F(t^{-m_\varphi},\varphi (t^{-m_\varphi},\xi_v (t^{-m_\varphi}))$ tends to a point $(a,0)\in \C^2$ as $t\mapsto 0$. This is impossible since $F$ is a diffeomorphism.

Now, let 
$$ 
J_\varphi:=a_\varphi P_\varphi \frac d{d\xi }Q_\varphi -b_\varphi Q_\varphi \frac 
d{d\xi }P_\varphi .$$ 
Taking  differentiation of $DF(t^{-m_\varphi},\varphi (t^{-m_\varphi},\xi )$, by (3) one can get that  
$$m_\varphi J(P,Q)t^{n_\varphi-2m_\varphi-1}     
=-J_\varphi t^{-a_\varphi-b_\varphi-1}+ \mbox { higher terms in }t. 
$$ 
Since $J(P,Q)\equiv const.\neq 0$,
$$ 
J_\varphi  \equiv \cases{-m_\varphi J(P,Q),&if $a_\varphi+b_\varphi+n_\varphi =2m_\varphi$ 
\cr 
        0,&if $a_\varphi+b_\varphi+n_\varphi >2m_\varphi.$\cr} 
$$  
If $J_\varphi  \equiv 0$, it must be that $P_\varphi^{-b_\varphi} =CQ_\varphi^{-a_\varphi}$ for $C\in\C^*$. This is impossible by Claim 1(b). Thus,
$J_\varphi =-m_\varphi J(P,Q)$. In particular, $P_\varphi$ and $Q_\varphi$ have only simple zeros. \endproof

\medskip

Now, we can complete the proof of the lemma.  By  Claim 2 the numbers $\alpha_u$ and $\beta_v$ are simple zero of $P_\varphi$ and $Q_\varphi$, respectively. Then, by Fact 2 there exists Newton-Puiseux expansions at infinity
$$\bar u(x^{\frac 1{m_\varphi}} )=\varphi (x,\alpha_u+\mbox{ lower terms in } x^{\frac 1{m_\varphi}}  ),$$
$$\bar v(x^{\frac 1{m_\varphi}} )=\varphi (x,\beta_v+\mbox{ lower terms in } x^{\frac 1{m_\varphi} } ),$$
for which $P(x,\bar u(x^{\frac1{m_\varphi}}))\equiv 0$,  $Q(x,\bar v(x^{\frac1{m_\varphi}}))\equiv 0$ and $\mult(\bar u)$  and  $\mult(\bar v )$ divide $ m_\varphi .$
Since $\mult(\bar u)=\deg P >\deg Q=\mult(\bar v)$ and 
$\alpha_u \neq \beta_v$, we get $\alpha_u\neq 0$, $\beta_v=0$ and $\deg P =m_\varphi$ . Hence,  $\deg Q \vert \deg P.$ \endproof

\medskip

\noindent {\bf 4. Proof of Automorphism Theorem.}
The proof uses Division Lemma and the following fact which is only an easy elementary excise on homogeneous polynomial. 

(*) {\it Let $f,g\in \C[x,y]$ be homogeneous. If $f_xg_y-f_yg_x\equiv 0$, then there is a homogeneous polynomial $h\in \C[x,y]$ with $\deg h=\gcd(\deg f, \deg g)$ such that }
$$ f=ah^{\frac {\deg f}{\deg h}}  \mbox{ and } g=ah^{\frac {\deg g}{\deg h}}, \ a,b\in \C^*.$$ (See, for example [E, Lemma 10.2.4, p 253]).

Given $F=(P,Q) \in GA(\C^2)$. Assume, for instance, $\deg P \geq \deg Q $ and $\deg P > 1$. Then, by the Division Lemma $\deg P =m \deg Q$,  and hence, by (*) $ \deg (P-cQ^m) <\deg P$ for a suitable number $c\in \C$. By induction one can find a finite sequence of automorphisms $\phi_i(x,y)$, $ i=1,\dots , k$ of the form $(x,y)\mapsto (x+cy^l,y)$ and $(x,y)\mapsto (x,y+cx^n)$ such that the components of  the map of $\phi_k\circ\phi_{k-1}\circ \dots \circ\phi_1\circ F$ are of degree $1$. Note that $\phi_i^{-1}$ has the form as those of $\phi_i$. Then, we get the automorphism Theorem. \endproof

\medskip 
 
\noindent {\bf Acknowlegments}: 
The author wishes to thank Prof. 
A.V. Essen and Prof. V.H. Ha  for their helps.  
 
\medskip 
 
\noindent {\bf References} 
 
\noindent [A]  S. S. Abhyankar,    {\it Expansion Techniques in Algebraic Geometry},    
Tata Institute of Fundamental Research, 1977.  

\noindent [AM] S. S. Abhyankar and T. T.  Moh,    {\it  Embeddings of the line in the 
plane},    J.  Reine Angew.  Math.  276 (1975),    148-166.  
 
\noindent [BCW] H.  Bass,    E.  Connell and D.  Wright,   {\it  The Jacobian conjecture: reduction of degree and formal expansion of the inverse},    Bull.  Amer.  Math.  Soc. (1982),    287-330.  
 
\noindent [BK] E.  Brieskorn,    H.  Knorrer,    {\it Ebene algebraische Kurven},    
Birkhauser,    Basel-Boston-Stuttgart 1981.  
    
\noindent [CK] J.Chadzy'nski and KT. Krasi'nski, {\it On a formula for the geometric 
degree and Jung' theorem }, Universitatis Iagellonicae Acta Mathematica, Fasciculus 
XXVIII, 1991, 81-84. 

\noindent [D] W. Dicks, {\it Automorphisms of the polynomial ring in two variables}, Publ. Sec. Math. Univ Autonoma Barcelona 27, (1983), 155-162.

\noindent [E] van den Essen,  Arno,  {\it Polynomial automorphisms and the Jacobian conjecture}.  (English.  English summary) Progress in Mathematics,  190.  Birkhauser Verlag,  Basel,  2000. 

\noindent [G] A. Gutwirth, {\it An inequality for certain pencils of plane curve}, Proc. Amer. Math. Soc. 12 (1961), 631-638.

\noindent [J] H. W. E.  Jung,   {\it  Uber ganze birationale Transformationen der 
Ebene },    J.  Reine Angew.  Math.  184 (1942),   161-174. 
  
\noindent [K]  O. Keller,     {\it Ganze Cremona-Transformationen},     Monatsh.     Mat.  
Phys.  47 (1939),     299-306.  
 
\noindent [Ku] W.  Van der Kulk,    {\it On polynomial rings in two variables},    Nieuw 
Arch.  Wisk.  (3) 1 (1953),    33-41.  

\noindent [MW] J. H.  McKay and S. S. Wang,   {\it  An elementary proof of the 
automorphism theorem for the polynomial ring in two variables },    J.  Pure Appl.  
Algebra 52 (1988),    91-102. 
 
\noindent [N] M. Nagata, {\it On Automorphism group of $k[x,y]$}, Lectures in Mathematics 5, Tokyo (1972).

\noindent [R] R. Rentschler, {\it Operations  du groupe additif sur le plan affine}, C.R.Acad. Sci. Paris Ser. A (1968), 384-387.

\noindent [Sh] I.R. Shafarevich, {\it On some infinite dimensional groups}, Rend.Math. Appl. 25 (1966), 208-212.

\bigskip

\noindent Hanoi Institute of Mathematics,
P.O. Box 631, Boho 10000,
Hanoi, Vietnam.\\
{\small E-mail: nvchau@thevinh.ntsc.ac.vn }

\end{document}